\documentclass[a4paper,3p,sort&compress,final]{elsarticle}
\usepackage[english]{babel}

\usepackage{hyperref}
\hypersetup{
  colorlinks=true,
  frenchlinks=false,
  pdfborder={0 0 0},
}

\usepackage{amssymb,amsmath,amsthm,amsfonts,dsfont, cleveref, graphicx}
\usepackage[dvipsnames]{xcolor}

\newtheorem{theorem}{Theorem}[section]
\newtheorem{proposition}[theorem]{Proposition}

\newtheorem{remark}[theorem]{Remark}

\newtheorem{conjecture}[theorem]{Conjecture}

\crefname{section}{section}{sections}
\crefname{subsection}{subsection}{subsections}
\Crefname{section}{Section}{Sections}
\Crefname{subsection}{Subsection}{Subsections}
\crefname{theorem}{Theorem}{Theorems}
\crefname{proposition}{Proposition}{Propositions}
\crefname{lemma}{Lemma}{Lemmata}
\crefname{definition}{Definition}{Definitions}
\crefname{conjecture}{Conjecture}{Conjectures}
\crefname{corollary}{Corollary}{Corollaries}
\Crefname{figure}{Figure}{Figures}

\crefformat{equation}{\textup{#2(#1)#3}}
\crefrangeformat{equation}{\textup{#3(#1)#4--#5(#2)#6}}
\crefmultiformat{equation}{\textup{#2(#1)#3}}{ and \textup{#2(#1)#3}}
{, \textup{#2(#1)#3}}{, and \textup{#2(#1)#3}}
\crefrangemultiformat{equation}{\textup{#3(#1)#4--#5(#2)#6}}%
{ and \textup{#3(#1)#4--#5(#2)#6}}{, \textup{#3(#1)#4--#5(#2)#6}}{, and \textup{#3(#1)#4--#5(#2)#6}}

\Crefformat{equation}{#2Equation~\textup{(#1)}#3}
\Crefrangeformat{equation}{Equations~\textup{#3(#1)#4--#5(#2)#6}}
\Crefmultiformat{equation}{Equations~\textup{#2(#1)#3}}{ and \textup{#2(#1)#3}}
{, \textup{#2(#1)#3}}{, and \textup{#2(#1)#3}}
\Crefrangemultiformat{equation}{Equations~\textup{#3(#1)#4--#5(#2)#6}}%
{ and \textup{#3(#1)#4--#5(#2)#6}}{, \textup{#3(#1)#4--#5(#2)#6}}{, and \textup{#3(#1)#4--#5(#2)#6}}

\crefformat{appendix}{#2#1#3}

\crefdefaultlabelformat{#2\textup{#1}#3}

\DeclareMathAlphabet{\mathpzc}{OT1}{pzc}{m}{it}

\newcommand{\norm}[1]{\left\|#1\right\|}
\newcommand{\vb}[1]{\mathbf{#1}}
\newcommand{\identity}{I}

\begin{document}

\author{Nick Dewaele}
\ead{nick.dewaele@kuleuven.be}
\address{KU Leuven, Department of Computer Science, Leuven, Belgium.}

\author{Paul Breiding\fnref{cor2}}
\ead{breiding@mis.mpg.de}
\address{Max-Planck-Institute for Mathematics in the Sciences, Leipzig, Germany.}
\fntext[cor2]{Funding: supported by the Deutsche Forschungsgemeinschaft (DFG) -- Projektnummer 445466444.}

\author{Nick Vannieuwenhoven\fnref{cor1}}
\ead{nick.vannieuwenhoven@kuleuven.be}
\address{KU Leuven, Department of Computer Science, Leuven, Belgium;\\ Leuven.AI, KU Leuven Institute for AI, B-3000 Leuven, Belgium.\vspace{-3em}}
\fntext[cor1]{Funding: supported by a Postdoctoral Fellowship of the Research Foundation---Flanders (FWO) with project 12E8119N.}

\title{Three decompositions of symmetric tensors have similar condition numbers}

\begin{abstract}
We relate the condition numbers of computing three decompositions of symmetric tensors: the canonical polyadic decomposition, the Waring decomposition, and a Tucker-compressed Waring decomposition. Based on this relation we can speed up the computation of these condition numbers by orders of magnitude.
\end{abstract}

\maketitle

\vspace{-1.75em}
\section{Introduction}

Many problems in machine learning and signal processing involve computing a decomposition of a \textit{symmetric tensor} \cite{Kolda2009,Anandkumar2014b}; an order-$D$ tensor $\mathpzc{A} = [a_{i_1,\dots,i_D}]_{i_1,\dots,i_D=1}^n \in \mathbb{R}^{n \times \dots \times n}$ is symmetric if its entries $a_{i_1,\dots,i_D}$
are invariant under all permutations of the indices $i_1,\dots,i_D$.
We establish a close connection between the numerical sensitivity of the following three increasingly structured decomposition problems associated with a symmetric tensor $\mathpzc{A}$:
\begin{enumerate}
 \item A \textit{canonical polyadic decomposition (CPD)} of $\mathpzc{A}$ expresses $\mathpzc{A}$ as a sum of $R$ (not necessarily symmetric) tensors of rank 1, where $R$ is minimal. In other words, $\mathpzc{A} = \sum_{r=1}^R \mathpzc{A}_r$ where $\mathpzc{A}_r = \alpha_r\ \vb{a}_r^{(1)}\otimes \cdots \otimes \vb{a}_r^{(D)}$,
 $\alpha_r \in \mathbb{R} \setminus \{0\}$ and each $\vb{a}_r^{(i)}$ is a point on the sphere $\mathbb{S}^{n - 1}=\{ \vb{a} \in \mathbb{R}^n \mid \|\vb{a}\|_2=1 \}$.
 \item A \textit{Waring decomposition (WD)} is a special case of the CPD where all summands are symmetric. That is, for $r=1,\dots,R$, we have that $\mathpzc{A}_r = \alpha_r\ \vb{a}_r^{\otimes D}$ where $\alpha_r \in \mathbb{R} \setminus \{0\}$, $\vb{a}_r\in \mathbb{S}^{n - 1}$, and
 $\vb{a}_r^{\otimes D}$ is the tensor product of $D$ copies of $\vb{a}_r$.
\item A \textit{$Q$-compressed Waring decomposition ($Q$-WD)} is defined as follows. A symmetric tensor $\mathpzc{A}$ can be represented in a minimal subspace by a symmetric Tucker decomposition \cite{Tucker1966,DeLathauwer2000}, i.e., $\mathpzc{A} = (Q,\dots,Q) \cdot \mathpzc{G}$ where $Q \in \mathbb{R}^{n \times m}$ has orthonormal columns and
$\mathpzc{G} \in \mathbb{R}^{m \times \dots \times m}$ is symmetric with $m < n$. We write this as $\mathpzc{A} = Q^{\otimes D} \mathpzc{G}$. If $\mathpzc{G}$ has a WD $\mathpzc{G} = \sum_{r=1}^R \mathpzc{G}_r$, then it can be converted to a WD $\mathpzc{A} = \sum_{r=1}^R Q^{\otimes D} \mathpzc{G}_r$.
\end{enumerate}
In all three cases, the summands are points on a smooth manifold $\mathcal{M}\subset\mathbb{R}^{n\times\dots\times n}$, so they are \textit{join decompositions} \cite{Breiding2018a}.
For the CPD, the summands lie on the \textit{Segre manifold} $\mathcal{S}_{n,D}$, for the WD they lie on \textit{Veronese manifold} $\mathcal{V}_{n, D}$  \cite{Landsberg2012}, and for the $Q$-WD they lie on the manifold $\mathcal{W}_{Q,D} = Q^{\otimes D}(\mathcal{V}_{m, D})$. In the remainder, we drop the subscripts on the manifolds if they are clear from the context.

We study the sensitivity of the summands in these three decompositions with respect to perturbations of~$\mathpzc{A}$. Consider a decomposition of $\mathpzc{A}$ with summands $\mathfrak{a} = (\mathpzc{A}_1,\dots,\mathpzc{A}_R) \in \mathcal{M}^{\times R}$, where $\mathcal{M}$ is one of the three manifolds described above and $\mathcal{M}^{\times R}$ is the product of $R$ copies of $\mathcal{M}$.
Under mild conditions \cite{Breiding2018a}, ${\mathfrak{a}}$ is an \emph{isolated} decomposition of $\mathpzc{A}$ and the addition map $\Sigma: \mathcal{M}^{\times R} \mapsto \mathbb{R}^{n\times\dots\times n},\; (\mathpzc{A}_1,\dots,\mathpzc{A}_R) \mapsto \mathpzc{A}_1+\dots+\mathpzc{A}_R$
admits a local inverse $\Sigma_{\mathfrak{a}}^{-1}$. In this case, the sensitivity of the decomposition with respect to $\mathpzc{A}$ can be measured by the condition number \cite{Rice1966}:
\begin{equation}
    \label{eq:defcond}
\kappa_{\mathcal{M}}(\mathpzc{A}_1,\dots,\mathpzc{A}_R) :=
\lim_{\delta \rightarrow 0} \;
\sup_{{\widetilde{\mathpzc{A}} \in \Sigma(\mathcal{M}^{\times R}),\, {\| \mathpzc{A} - \widetilde{\mathpzc{A}} \|} \le \delta}} \frac{ \| \Sigma^{-1}_{\mathfrak{a}}(\mathpzc{A}) - \Sigma^{-1}_{\mathfrak{a}} (\widetilde{\mathpzc{A}}) \| }{ \|\mathpzc{A} - \widetilde{\mathpzc{A}} \|},
\end{equation}
where $\norm{\cdot}$ is the Euclidean or Frobenius norm. If $\mathfrak{a}$ is not isolated, $\kappa_{\mathcal{M}}(\mathpzc{A}_1,\dots,\mathpzc{A}_R) := \infty$.

Suppose $\mathpzc{A}$ has a $Q$-WD $\mathfrak{a} = (\mathpzc{A}_1,\dots,\mathpzc{A}_R)$. It can also be regarded as a WD or CPD by ignoring symmetry or subspace constraints. We investigate the relationship between the condition numbers of these three problems at $\mathfrak{a}$. Since $\mathcal{W} \subseteq \mathcal{V} \subseteq \mathcal{S}$, it follows from \cref{eq:defcond} that $\kappa_{\mathcal{W}}(\mathfrak{a}) \le \kappa_{\mathcal{V}}(\mathfrak{a}) \le \kappa_{\mathcal{S}}(\mathfrak{a})$.
Similarly to recent findings \cite{Dewaele2021},
we show the following results for the WD:

\begin{theorem}
    \label{prop:WaringCondInvariance}
    Let $\mathpzc{G} = \mathpzc{G}_1 + \dots + \mathpzc{G}_R \in \mathbb{R}^{m \times \dots \times m}$ be a WD of an order-$D$ tensor.
    \begin{enumerate}
    \item Take $Q \in \mathbb{R}^{n \times m}$ with orthonormal columns and set $\mathpzc{A}_r = Q^{\otimes D} \mathpzc{G}_r$, for $r=1,\dots,R$. Then
    $$
    \kappa_{\mathcal{W}_{Q,D}}(\mathpzc{A}_1,\dots,\mathpzc{A}_R) \le
    \kappa_{\mathcal{V}_{n,D}}(\mathpzc{A}_1,\dots,\mathpzc{A}_R)
    \le
    \sqrt{D} \cdot
    \kappa_{\mathcal{V}_{m,D}}(\mathpzc{G}_1,\dots,\mathpzc{G}_R) = \sqrt{D} \cdot \kappa_{\mathcal{W}_{Q,D}}(\mathpzc{A}_1,\dots,\mathpzc{A}_R)
    .$$
    \item Let $U \in \mathbb{R}^{\ell \times m}$ have orthonormal columns and $\mathpzc{B}_r := U^{\otimes D} \mathpzc{G}_r$ for all $r$. If $\min(\ell, n) > m$, then
    $
    \kappa_{\mathcal{V}_{n,D}}(\mathpzc{A}_1,\dots,\mathpzc{A}_R)
    =
    \kappa_{\mathcal{V}_{\ell,D}}(\mathpzc{B}_1,\dots,\mathpzc{B}_R)
    ;$
i.e., the condition number is invariant under non-minimal symmetric Tucker compressions.
\end{enumerate}
\end{theorem}

Numerical evidence indicates a stronger connection, which can be proved in the rank-2 case:

\begin{conjecture}
    \label{conjecture:WaringCondition}
    If $\mathpzc{A} = \sum_{r=1}^R \mathpzc{A}_r$ is a WD of an order-$D$ symmetric tensor $\mathpzc{A} \in \mathbb{R}^{n\times\dots\times n}$ with $D \geq 3$, then
    $
    \kappa_{\mathcal{V}}(\mathpzc{A}_1,\dots,\mathpzc{A}_R)
    =
    \kappa_{\mathcal{S}}(\mathpzc{A}_1,\dots,\mathpzc{A}_R)
    .$
\end{conjecture}

\begin{proposition}
    \label{prop:waringCondSpecialCase}
\Cref{conjecture:WaringCondition} holds for $R \le 2$.
\end{proposition}

In conjunction with \cite[Theorem 5.1]{Dewaele2021}, \cref{conjecture:WaringCondition} would imply that $\kappa_{\mathcal{W}}(\mathfrak{a}) = \kappa_{\mathcal{V}}(\mathfrak{a}) = \kappa_{\mathcal{S}}(\mathfrak{a})$ for any $Q$-WD $\mathfrak{a}$,
which is sharper than \cref{prop:WaringCondInvariance}. This entails that the supremum in \cref{eq:defcond} applied to the Segre manifold (i.e., $\mathcal{M} = \mathcal{S}$) can be attained locally with a perturbation $\widetilde{\mathpzc{A}} \in \Sigma(\mathcal{W}^{\times R})$.

A practical consequence of \cref{prop:WaringCondInvariance} relates to the following procedure from \cite{Breiding2018a} to compute the condition number. Let $\mathcal{M}$ be either $\mathcal{S}$ or $\mathcal{V}$.
Let the matrix $T_{\mathpzc{A}_r}^{\mathcal{M}}$ contain as columns an orthonormal basis of $T_{\mathpzc{A}_r} \mathcal{M}$ for $r=1,\dots,R$. Then the condition number is characterized by the \textit{Terracini matrix} $T_{\mathpzc{A}_1,\dots,\mathpzc{A}_R}^{\mathcal{M}}$:
\begin{equation}
    \label{eq:GeneralTerraciniMtx}
T_{\mathpzc{A}_1,\dots,\mathpzc{A}_R}^{\mathcal{M}} := \left[ T_{\mathpzc{A}_1}^\mathcal{M} \quad \cdots \quad T_{\mathpzc{A}_R}^{\mathcal{M}} \right]
\quad\text{and}\quad
\kappa_{\mathcal{M}}(\mathpzc{A}_1,\dots,\mathpzc{A}_R) = \sigma_{\min}(T_{\mathpzc{A}_1,\dots,\mathpzc{A}_R}^{\mathcal{M}})^{-1},
\end{equation}
where $\sigma_{\min}(A)$ is the smallest singular value of $A$.
Consider a WD $\mathpzc{A} = \sum_{r=1}^R \mathpzc{A}_r$ with $\mathpzc{A}_r = \alpha_r \vb{a}_r^{\otimes D}$
for some $\alpha_r \in \mathbb{R} \setminus \{0\}$ and $\vb{a}_r \in \mathbb{S}^{n - 1}$. For this decomposition the Terracini matrices for the CPD and the WD, respectively, are given as follows: for any two matrices $X$ and $A$, let $X \otimes_d A^{\otimes D - 1} := A^{\otimes d - 1} \otimes X \otimes A^{\otimes D - d}$.
Let $U(\vb{a}_r)$ be any orthonormal basis of $T_{\vb{a}_r} \mathbb{S}^{n - 1} = \vb{a}_r^\perp$. Then the Terracini matrices are
\begin{equation}
    \label{eq:TerraciniCPDandWaring}
T^{\mathcal{S}}_{\mathpzc{A}_1,\dots,\mathpzc{A}_R} =
\left[
\vb{a}_r^{\otimes D} \quad
\Big[
    U_r \otimes_d  \vb{a}_r^{\otimes D - 1}
\right]_{d=1}^D
\Big]_{r=1}^R
\quad\text{and}\quad
T^{\mathcal{V}}_{\mathpzc{A}_1,\dots,\mathpzc{A}_R} =
\Big[
\vb{a}_r^{\otimes D} \quad
\frac{1}{\sqrt{D}} \sum_{d=1}^D U(\vb{a}_r) \otimes_d  \vb{a}_r^{\otimes D - 1}
\Big]_{r=1}^R.
\end{equation}
A major implication of \cref{prop:WaringCondInvariance} is that we can speed up the computation of $\kappa_{\mathcal{V}_{n,D}}(\mathpzc{A}_1,\dots,\mathpzc{A}_R)$. 
Assuming $n > R$ and $\mathpzc{A}_r = \alpha_r \vb{a}_r^{\otimes D}$ with $\alpha_r \ne 0$ and $\vb{a}_r \in \mathbb{S}^{n - 1}$, the following computes $\kappa_{\mathcal{V}_{n,D}}(\mathpzc{A}_1,\dots,\mathpzc{A}_R)$
\begin{enumerate}
    \item Compute a thin singular value decomposition $[\vb{a}_r]_{r=1}^R = Q \Sigma V^T$ and set $[\vb{g}_r]_{r=1}^R := \Sigma V^T \in \mathbb{R}^{m \times R}$.
    \item Construct $\vb{b}_r = [\vb{g}_r^T \quad 0]^T \in \mathbb{R}^{\ell}$ where $\ell = m + 1$ and set $\mathpzc{B}_r = \alpha_r \vb{b}_r^{\otimes D}$ for each $r$.
    \item Construct $T_{\mathpzc{B}_1,\dots,\mathpzc{B}_R}^{\mathcal{V}}$ as in \cref{eq:TerraciniCPDandWaring} and compute $\kappa_{\mathcal{V}}(\mathpzc{B}_1,\dots,\mathpzc{B}_R)$ by appling \cref{eq:GeneralTerraciniMtx}.
\end{enumerate}
Steps 1-2 give one possible choice of $Q$ and $U = [\identity \quad 0]^T$ and $\mathpzc{G}_r = \alpha_r \vb{g}_r^{\otimes D}$ that satisfy \cref{prop:WaringCondInvariance}.
A Julia \cite{Bezanson2017} implementation of this method is provided along with the arXiv version of this manuscript.
Since $T^{\mathcal{V}}_{\mathpzc{B}_1,\dots,\mathpzc{B}_R} \in \mathbb{R}^{\ell^D \times R\ell}$ and $\ell \le R + 1$, step 3 can be performed in $\mathcal{O}(R^{D + 4})$ operations, adding to the $\mathcal{O}(nR^2)$ cost of step 1. Applying \cref{eq:GeneralTerraciniMtx} to $T^{\mathcal{V}}_{\mathpzc{A}_1,\dots,\mathpzc{A}_R} \in \mathbb{R}^{n^D \times Rn}$ would involve $\mathcal{O}(n^{D+2}R^2)$ operations. The algorithm can reach significant speedups if $n \gg R$. For instance, we applied the Julia code to a WD with $(n,D,R,\ell) = (100, 3, 10, 11)$ on an Intel Xeon CPU E5-2697 v3 running on 8 cores and 126GB memory. The computation times were $115.6$ and $0.0092$ seconds for the original and improved algorithm, respectively.

\section{Condition number of a $Q$-WD}
\label{sec:condCompressedWD}

In this section, we prove \cref{prop:WaringCondInvariance} based on the following insight: $\Sigma(\mathcal{V}^{\times R})$ is locally a manifold whose tangent space is decomposed as $\mathbb{T} \oplus \mathbb{T}^{\perp}$ where $\mathbb{T}$ is the tangent space to $\Sigma(\mathcal{W}^{\times R})$ and $\mathbb{T}^{\perp}$ is its orthogonal complement. As long as $n > m$, the effect of the worst perturbation to $\mathpzc{A}$ inside $\mathbb{T}^{\perp}$ is independent of $n$ and can be bounded as in the first statement. From this, the second statement follows as well.

\begin{proof}[Proof of \cref{prop:WaringCondInvariance}]
The first inequality follows from the inclusion $\mathcal{W}_{Q,D} \subseteq \mathcal{V}_{n,D}$. The last follows from the fact that $Q^{\otimes D}$ is an isometry between $\mathcal{V}_{m,D}$ and $\mathcal{W}_{Q,D}$. It remains to show the middle inequality. If $m = n$, $Q$ is an orthogonal change of basis, which preserves the condition number. Thus, we assume $n > m$.

For each $r$, let $\mathpzc{G}_r = \alpha_r \vb{g}_r^{\otimes D}$ with $\alpha_r \in \mathbb{R} \setminus \{0\}$ and $\vb{g}_r \in \mathbb{S}^{m - 1}$, let $\vb{a}_r = Q \vb{g}_r$ and define $U_r$ so that the matrix~$[\vb{g}_r \quad U_r] \in \mathbb{R}^{m \times m}$ is orthogonal.
Construct $T_{\mathpzc{G}_r}^{\mathcal{V}}$ by applying \cref{eq:TerraciniCPDandWaring} to $\mathpzc{G}_r$.
Complete $Q$ to an orthonormal basis $[Q \quad Q_\perp]$ of $\mathbb{R}^{n}$.
The columns of $U(\vb{a}_r) := [Q U_r \quad Q_\perp]$ form an orthonormal basis of $T_{\vb{a}_r} \mathbb{S}^{n - 1}$. Substituting this into \cref{eq:TerraciniCPDandWaring} gives
\begin{align*}
T_{\mathpzc{A}_1,\dots,\mathpzc{A}_R}^{\mathcal{V}_{n,D}}
= \begin{bmatrix} T_r & T_r^\perp \end{bmatrix}_{r=1}^R
\text{ with }
T_r
=
\begin{bmatrix}
    \vb{a}_r^{\otimes D}
&
    \frac{1}{\sqrt{D}}
    \left(\sum_{d=1}^{D} Q U_r \otimes_d \vb{a}_r^{\otimes D - 1}\right)
\end{bmatrix}\text{ and }
T_r^\perp
=
    \frac{1}{\sqrt{D}}
    \sum_{d=1}^{D}
    Q_\perp \otimes_d \vb{a}_r^{\otimes D - 1}
.\end{align*}

Since $\vb{a}_r = Q \vb{g}_r$, we have $T_r = Q^{\otimes D} T_{\mathpzc{G}_r}^{\mathcal{V}_{m, D}}$.
Thus, up to a column permutation, $T_{\mathpzc{A}_1,\dots,\mathpzc{A}_R}^{\mathcal{V}_{n,D}}$ is the horizontal concatenation of $Q^{\otimes D} T_{\mathpzc{G}_1,\dots,\mathpzc{G}_R}^{\mathcal{V}_{m, D}}$ and $T^\perp := \left[\begin{smallmatrix} T_1^\perp & \dots & T_R^\perp \end{smallmatrix}\right]$.
The column spaces of $Q^{\otimes D}$ and $T^\perp$ are orthogonal, so that the singular values of $T_{\mathpzc{A}_1,\dots,\mathpzc{A}_R}^{\mathcal{V}_{n,D}}$ are the union of those of $Q^{\otimes D} T_{\mathpzc{G}_1,\dots,\mathpzc{G}_R}^{\mathcal{V}_{m, D}}$ and $T^\perp$ separately.
Since $Q$ has orthonormal columns, $Q^{\otimes D} T_{\mathpzc{G}_1,\dots,\mathpzc{G}_R}^{\mathcal{V}_{m, D}}$ has the same singular values as $T_{\mathpzc{G}_1,\dots,\mathpzc{G}_R}^{\mathcal{V}_{m, D}}$, so it suffices to show $\sigma_{\min}(T^\perp) \ge \sigma_{\min}(T_{\mathpzc{G}_1,\dots,\mathpzc{G}_R}^{\mathcal{V}_{m,D}}) / \sqrt{D}$.

To do this, we compute $(T^\perp)^T T^\perp = [(T_{r_1}^\perp)^T (T_{r_2}^\perp)]_{r_1,r_2=1}^R$, where the block at $(r_1,r_2)$ is
\begin{equation}
    \label{eq:gramianOrthComplement}
\frac{1}{D}
\left(\sum_{d=1}^{D} Q_\perp \otimes_d \vb{a}_{r_1}^{\otimes D - 1} \right)^T
\left( \sum_{d=1}^{D} Q_\perp \otimes_d \vb{a}_{r_2}^{\otimes D - 1} \right)
=
\left\langle
\vb{a}_{r_1},\,
\vb{a}_{r_2}
\right\rangle^{D - 1} \identity_{n - m}
=
\left\langle
    \vb{g}_{r_1}
,\,
    \vb{g}_{r_2}
\right\rangle^{D - 1} \identity_{n - m}
.\end{equation}
Consider two modifications of $T^\perp$ that preserve the singular values: first, let
$\hat{T}^\perp := [\identity_{n - m} \otimes \vb{g}_r^{\otimes D - 1}]_{r=1}^R$, then by \cref{eq:gramianOrthComplement}, we have $(\hat{T}^\perp)^T\hat{T}^\perp = (T^\perp)^T (T^\perp)$, so they have the same singular values.
Second, if we define $\widetilde{T}^\perp := [[\vb{g}_r \quad U_r] \otimes \vb{g}_r^{\otimes D - 1}]_{r=1}^R$,
then $\widetilde{T}^\perp$ and $\hat{T}^\perp$ also have the same singular values, since $[\vb{g}_r \quad U_r]$ and $\identity_{n - m}$ are orthogonal up to multiplicities \cite[lemma 5.3]{Dewaele2021}.
Hence, we can proceed with $\widetilde{T}^\perp$ instead of $T^\perp$.
Similarly, we modify $T_{\mathpzc{G}_1,\dots,\mathpzc{G}_R}^{\mathcal{V}}$. Scaling up all its columns of the form $\vb{g}_r^{\otimes D}$ by $\sqrt{D}$ gives
$$
\widetilde{T}_{\mathpzc{G}_1,\dots,\mathpzc{G}_R}^{\mathcal{V}} :=
\left[
    \sqrt{D} \vb{g}_r^{\otimes D}
    \quad
    \frac{1}{\sqrt{D}}
    \sum_{d=1}^D U_r \otimes_d \vb{g}_r^{\otimes D - 1}
\right]_{r=1}^R
=
\left[
    \frac{1}{\sqrt{D}} \sum_{d=1}^D [\vb{g}_r \quad U_r] \otimes_d \vb{g}_r^{\otimes D - 1}
\right]_{r=1}^R
,$$
i.e., $\widetilde{T}_{\mathpzc{G}_1,\dots,\mathpzc{G}_R}^{\mathcal{V}} = T_{\mathpzc{G}_1,\dots,\mathpzc{G}_R}^{\mathcal{V}} \Delta$ where $\Delta$ is diagonal and $\sigma_{\min}(\Delta) = 1$. Hence, $\sigma_{\min}(\widetilde{T}_{\mathpzc{G}_1,\dots,\mathpzc{G}_R}^{\mathcal{V}}) \ge \sigma_{\min}(T_{\mathpzc{G}_1,\dots,\mathpzc{G}_R}^{\mathcal{V}})$.

To compare the singular values of $\widetilde{T}^\perp$ and $\widetilde{T}_{\mathpzc{G}_1,\dots,\mathpzc{G}_R}^{\mathcal{V}}$, take the singular vector $\vb{v} = [\vb{v}_r \in \mathbb{R}^{m}]_{r=1}^R$ of $\widetilde{T}^\perp$ corresponding to the smallest singular value and compute
$$
\widetilde{T}_{\mathpzc{G}_1,\dots,\mathpzc{G}_R}^{\mathcal{V}_{m,D}} \vb{v}
=
\sum_{r=1}^R \left(
    \frac{1}{\sqrt{D}} \sum_{d=1}^D [\vb{g}_r \quad U_r] \otimes_d \vb{g}_r^{\otimes D - 1}
\right) \vb{v}_r
=
\frac{1}{\sqrt{D}}
\sum_{d=1}^D \sum_{r=1}^R
     ([\vb{g}_r \quad U_r] \vb{v}_r) \otimes_d \vb{g}_r^{\otimes D - 1}
.$$
Since all the summands in the outer sum have the same norm, the triangle inequality gives
$$
\left\| \widetilde{T}_{\mathpzc{G}_1,\dots,\mathpzc{G}_R}^{\mathcal{V}_{m,D}} \vb{v} \right\|
\le
\sqrt{D}
\left\|
    \sum_{r=1}^R ([\vb{g}_r \quad U_r] \vb{v}_r) \otimes \vb{g}_r^{\otimes D - 1}
\right\|
=
\sqrt{D}\,
\|
    \left[\left[ \vb{g}_r \quad U_r \right] \otimes \vb{g}_r^{\otimes D - 1} \right]_{r=1}^R \vb{v}
\|
= \sqrt{D} \cdot \sigma_{\min}(\widetilde{T}^\perp).
$$
As $\sigma_{\min}(T_{\mathpzc{G}_1,\dots,\mathpzc{G}_R}^{\mathcal{V}}) \le \sigma_{\min} (\widetilde{T}_{\mathpzc{G}_1,\dots,\mathpzc{G}_R}^{\mathcal{V}}) \le \bigl\| \widetilde{T}_{\mathpzc{G}_1,\dots,\mathpzc{G}_R}^{\mathcal{V}} \vb{v} \bigr\|$, and $\sigma_{\min}(\widetilde{T}^\perp) = \sigma_{\min}(T^\perp)$
this gives the desired bound.

For the second statement, recall that the singular values of $T_{\mathpzc{A}_1,\dots,\mathpzc{A}_R}^{\mathcal{V}_{n,D}}$ are the union of those of $T_{\mathpzc{G}_1,\dots,\mathpzc{G}_R}^{\mathcal{V}_{m,D}}$ and those of $\widetilde{T}^\perp$ whenever $n > m$. Observe that both of these matrices are independent of $n$ and $Q$. Hence, applying the above calculation to $\mathpzc{B}_1,\dots,\mathpzc{B}_R\in \mathcal{V}_{\ell, D}\subset \mathbb R^{\ell \times \cdots \times \ell}$ and orthogonal $U\in\mathbb R^{\ell\times m}$ under the assumption $\ell>m$
would reveal the same singular values.
\end{proof}

\section{Equivalence between the CPD and WD}

\Cref{conjecture:WaringCondition} is a stronger statement than \cref{prop:WaringCondInvariance}, but it seems too challenging to show in general. We present a proof for the case where $R=2$ and present numerical evidence for the general case.

\begin{proof}[Proof of \cref{prop:waringCondSpecialCase}]
For $R = 1$, both condition numbers are equal to $1$ by \cref{eq:GeneralTerraciniMtx,eq:TerraciniCPDandWaring}.
For $R=2$, the proof comprises computing the singular values of \cref{eq:GeneralTerraciniMtx} for the CPD.
Let $\mathpzc{A}_1 = \lambda_1 \vb{u}^{\otimes D}$ and $\mathpzc{A}_2 = \lambda_2 \vb{v}^{\otimes D}$ with $\vb{u}, \vb{v} \in \mathbb{S}^{n - 1}$ and $\lambda_1,\lambda_2 \ne 0$.
Let $U, V \in \mathbb{R}^{n \times (n - 1)}$ be orthonormal bases of $T_{\vb{u}}\mathbb{S}^{n - 1} = \vb{u}^\perp$ and $T_{\vb{v}} \mathbb{S}^{n - 1} = \vb{v}^\perp$, respectively.
Applying \cref{eq:TerraciniCPDandWaring} and using as before the notation $X \otimes_d A^{\otimes D - 1} := A^{\otimes d - 1} \otimes X \otimes A^{\otimes D - d}$ gives
$$
T_{\mathpzc{A}_1,\mathpzc{A}_2}^{\mathcal{S}} = \left[
    \vb{u}^{\otimes D}
    \quad
    \left[
        U \otimes_d \vb{u}^{\otimes D - 1}
    \right]_{d=1}^D
    \quad
    \vb{v}^{\otimes D}
    \quad
    \left[
        V \otimes_d \vb{v}^{\otimes D - 1}
    \right]_{d=1}^D
\right]
.$$
Next, define the vectors $
\vb{q}_{D}^j := \bigl[
    \tfrac{1}{\sqrt{j (j+1)}} \vb{1}_j^T
    \quad \tfrac{-j}{\sqrt{j (j+1)}} \quad
    \vb{0}_{D - j - 1}^T
\bigr]^T
\in \mathbb{R}^{D}
$
where $\vb{1}_N, \vb{0}_N\in\mathbb R^N$ are the vectors consisting of ones and zeros, respectively. We set
$
Q_D := \left[
    \frac{1}{\sqrt{D}} \vb{1}_D \,\,\quad \vb{q}_D^1 \quad\dots\quad \vb{q}_D^{D-1}
\right]\in\mathbb R^{D\times D}
.$ This matrix is called \emph{Helmert's orthogonal matrix} \cite{Helmert1876}; the rows of its right $D\times (D-1)$ submatrix are the vertices of a regular simplex in $\mathbb R^{D-1}$.
We transform $T^{\mathcal{S}}_{\mathpzc{A}_1,\mathpzc{A}_2}$ into a matrix $\widetilde{T}^{\mathcal{S}}$ with the same singular values using an orthogonal change of basis:  $\widetilde{T}^{\mathcal{S}} :=
T^{\mathcal{S}}_{\mathpzc{A}_1,\mathpzc{A}_2} \mathrm{diag}(1, \identity \otimes Q_D, 1, \identity \otimes Q_D)$. This gives
\begin{align*}
&\widetilde{T}^{\mathcal{S}} = \left[
    \vb{u}^{\otimes D} \quad
    S_{\vb{u}} \quad S_{\vb{u}, \perp}^1 \quad \dots \quad S_{\vb{u}, \perp}^{D-1}
    \quad
    \vb{v}^{\otimes D} \quad
    S_{\vb{v}} \quad S_{\vb{v}, \perp}^1 \quad \dots \quad S_{\vb{v}, \perp}^{D-1}
\right]\\
\text{where} \quad
&S_{\vb{u}} = \frac{1}{\sqrt{D}} \sum_{d=1}^D(U \otimes_d \vb{u}^{\otimes D - 1})
\quad\text{and}\quad
\quad S_{\vb{u}, \perp}^j = \sum_{i=1}^{j+1} (\vb{q}_D^j)_i (U \otimes_i \vb{u}^{\otimes D - 1})
\end{align*}
and analogously for $\vb{v}$. After rearranging the blocks, we get the following partition of $\widetilde{T}^{\mathcal{S}}$:
\begin{equation}
    \label{eq:TerraciniSegreBlocks}
\widetilde{T}^{\mathcal{S}} \cong
\begin{bmatrix} T^{\mathcal{V}} & T^1_\perp & \dots & T^{D-1}_\perp
\end{bmatrix}
\,\,\text{ where }\,\,
T^{\mathcal{V}}_{\mathpzc{A}_1,\mathpzc{A}_2} = \begin{bmatrix}
    \vb{u}^{\otimes D} & S_{\vb{u}} &
    \vb{v}^{\otimes D} & S_{\vb{v}}
\end{bmatrix}
\,\,\text{ and }\,\,
T^d_\perp :=
\begin{bmatrix}
    S_{\vb{u},\perp}^d & S_{\vb{v},\perp}^d
\end{bmatrix}
\end{equation}
in which we recognise \cref{eq:TerraciniCPDandWaring}.
Now, we will show that these $D$ blocks are pairwise orthogonal, so that the singular values of $\widetilde{T}^{\mathcal{S}}$ are the union of the singular values of the blocks.
To see this, we compute $(S_{\vb{u},\perp}^j)^T (S_{\vb{v},\perp}^{j^\prime})$. Let~$\alpha := \left\langle \vb{u},\,\vb{v} \right\rangle$. Without loss of generality, assume $j < j'$. Let
$$x_i = (U \otimes_i \vb{u}^{\otimes D - 1})\quad \text{and}\quad y_{i'} = (V \otimes_{i'} \vb{v}^{\otimes D - 1});\quad  \beta_{=} := \alpha^{D - 1} U^T V \quad\text{and}\quad \beta_{\neq} := \alpha^{D - 2} U^T \vb{v} \vb{u}^T V.$$
Note that $x_i^T y_{i'}  = \beta_{=}$, if $i = i'$, and that $x_i^T y_{i'}  = \beta_{\neq}$, if $i \neq i'$.
Up to the constant $(j(j+1)j'(j'+1))^{-\frac{1}{2}}$, the inner products $(S_{\vb{u},\perp}^j)^T (S_{\vb{v},\perp}^{j^\prime})$ are
\begin{equation}
\label{eq:symComplementInnerProduct}
\left(
    x_1 + \dots + x_j - jx_{j+1}
\right)^T \left(
    y_1 + \dots + y_j + \dots + y_{j'} - j' y_{j' + 1}
\right),
\end{equation}
which is a linear form in $\beta_{=}$ and $\beta_{\ne}$.
First, we calculate the terms in \cref{eq:symComplementInnerProduct} involving the case $i = i^\prime$.
There are $j$ terms of the form $x_i^T y_i$ with $i \le j < j'$ and one of the form $-jx_{j+1} y_{j+1}$. Adding these terms together shows that the coefficient of $\beta_{=}$ is zero.
Second, we identify all terms in \cref{eq:symComplementInnerProduct} where the coefficient of $\beta_{\ne}$ is positive. For each $x_i$ with $i \le j$, there are $j' - 1$ terms $x_i^T y_{i'}$ with $i' \ne i$. One more term has a positive coefficient of $\beta_{\ne}$, which is $jj' x_{j+1} y_{j'+1}$. Together, these terms add up to $j(j'-1) \beta_{\neq} + jj' \beta_{\neq}$.
Third, we accumulate the negative coefficients of $\beta_{\neq}$, which involve either $x_{j+1}$ or $y_{j'+1}$. For $x_{j+1}$, there are $j' - 1$ terms $y_{i'}$ with $j+1 \ne i' \le j'$. For $y_{j' + 1}$, there are $j$ terms $x_i$ with $i \le j$. Hence, the terms with a negative coefficient of $\beta_{\neq}$ add up to $-j(j'-1) \beta_{\neq} - jj' \beta_{\neq}$. This means the terms involving $\beta_{\neq}$ also vanish. Therefore, all inner products $(S_{\vb{u},\perp}^j)^T (S_{\vb{v},\perp}^{j^\prime})$ vanish for $j \ne j'$.

Furthermore, the columns of $T^{\mathcal{V}}_{\mathpzc{A}_1,\mathpzc{A}_2}$ are symmetric tensors. The space of symmetric tensors is the linear span of the Veronese manifold $\mathcal V := \{\alpha\vb{x}^{\otimes D} \mid \alpha\in\mathbb R\setminus\{0\}, \vb{x} \in \mathbb{S}^{n-1} \}$.
Since $\sum_{i=1}^{j+1} (\vb{q}_D^j)_i = 0$, we have $(\vb{x}^{\otimes D})^T S_{\vb{u},\perp}^j = \sum_{i=1}^{j+1} (\vb{x}^T \vb{u})^{D - 1} \vb{x}^T U (\vb{q}_D^j)_i = 0$, so that the columns of $S_{\vb{u},\perp}^j$ and $T^{\mathcal{V}}_{\mathpzc{A}_1,\mathpzc{A}_2}$ are pairwise orthogonal. We can therefore conclude that \cref{eq:TerraciniSegreBlocks} partitions $\widetilde{T}^{\mathcal{S}}$ into pairwise orthogonal blocks.

Next, we compute all singular values of $\widetilde{T}^{\mathcal{S}}$ by computing the singular values of the blocks in \cref{eq:TerraciniSegreBlocks} separately.
Using the same notation as before, we compute the blocks of $(T_\perp^j)^TT_\perp^j$:
$$
(S_{\vb{u},\perp}^j)^T (S_{\vb{v},\perp}^j)
= \frac{1}{j(j+1)} \left(
    x_1 + \dots + x_j - jx_{j+1}
\right)^T \left(
    y_1 + \dots + y_j - jy_{j+1}
\right)
= \frac{1}{j(j+1)} \left(
    a + b + c + d
\right)
$$
where
$$
\begin{bmatrix}
    a \\
    b \\
    c \\
    d \\
\end{bmatrix}
=
\begin{bmatrix}
    (x_1 + \dots + x_j)^T (y_1 + \dots + y_j) \\
    - (x_1 + \dots + x_j)^T jy_{j+1} \\
    - j x_{j+1}^T (y_1 + \dots + y_j) \\
    j^2 x_{j+1}^T y_{j+1}
\end{bmatrix}
=
\begin{bmatrix}
    j \beta_{=} + (j^2 - j) \beta_{\neq} \\
    -j^2 \beta_{\neq} \\
    -j^2 \beta_{\neq} \\
    j^2 \beta_{=}
\end{bmatrix}
.$$
This gives $
(S_{\vb{u},\perp}^j)^T (S_{\vb{v},\perp}^j)
=
\beta_{=} - \beta_{\neq} = \alpha^{D-1} U^T V - \alpha^{D-2} U^T \vb{v}\vb{u}^T V
$. Hence, the Gramian of $T_\perp^j$ is
$$
G_\perp :=
(T_\perp^j)^T T_\perp^j
= \begin{bmatrix}
    \identity_{n-1}
    & \alpha^{D-1} U^T V - \alpha^{D-2} U^T \vb{v}\vb{u}^T V \\
    \alpha^{D-1} V^T U - \alpha^{D-2} V^T \vb{u}\vb{v}^T U
    & \identity_{n-1}
\end{bmatrix},
$$
which this is independent of $j$.
The Gramian of $T^{\mathcal{V}}_{\mathpzc{A}_1,\mathpzc{A}_2}$ is
$$
G_S :=
(T^{\mathcal{V}}_{\mathpzc{A}_1,\mathpzc{A}_2})^T
T^{\mathcal{V}}_{\mathpzc{A}_1,\mathpzc{A}_2}
=
   \begin{bmatrix}
    1 & 0 & \alpha^D & \sqrt{D} \alpha^{D-1} u^T V \\
    \times & \identity_{n-1} & \sqrt{D} \alpha^{D-1} U^T v  & \alpha^{D-1} U^T V + (D-1) \alpha^{D-2} U^T v u^T V \\
    \times  & \times & 1  & 0 \\
    \times & \times & \times & \identity_{n-1} \\
   \end{bmatrix},
$$
where each $\times$ should be replaced by the transpose of corresponding element in the upper diagonal part.

To continue, we exploit the liberty of choosing the bases $U$ and $V$ for the orthogonal complements $u^\perp$ and $v^\perp$ respectively. By planar geometry, we can choose these bases such that $U \vb{e}_1 = \frac{v - \alpha u}{\|v - \alpha u\|}$, $V \vb{e}_1 = \frac{u - \alpha v}{\|u - \alpha v\|}$ and $U \vb{e}_j = V \vb{e}_j$ for all $j=2,\dots,n-1$. Consequently,
$
 U^T \vb{v} = \sqrt{1 - \alpha^2} \vb{e}_1
$,
$V^T \vb{u} = \sqrt{1-\alpha^2} \vb{e}_1$, and $U^T V = \mathrm{diag}( -\alpha, 1, \dots, 1).
$
Plugging these into $G_\perp$, we get
\begin{align*}
 G_\perp &=
 \begin{bmatrix}
  \identity_{n-1}           & \mathrm{diag}(-\alpha^D - \alpha^{D-2}(1 - \alpha^2),\alpha^{D-1},\dots,\alpha^{D-1}) \\
  \times & \identity_{n-1}
 \end{bmatrix}
 = \identity_{2(n-1)} +
 \begin{bmatrix}
    0 & A_\perp \\
    A_\perp & 0
 \end{bmatrix},
\end{align*}
where $A_\perp := \mathrm{diag}(-\alpha^{D-2}, \alpha^{D-1}, \dots, \alpha^{D-1})$. Recall that the eigenvalues of $\left[\begin{smallmatrix} 0 & A^T \\ A & 0 \end{smallmatrix}\right]$ are $\pm \sigma(A)$, where $\sigma$ are the singular values of $A$. Therefore, the eigenvalues of $G_\perp$ are
$\lambda(G_\perp) = \{ 1 \pm \alpha^{D-1}, 1 \pm \alpha^{D-2} \}$
. We only need the extreme eigenvalues, which are $1 \pm \alpha^{D-2}$ since $| \alpha | \le 1$. For $G_S$, we obtain
\begin{align*}
G_S = \begin{bmatrix}
    1 & 0 & \alpha^D & \sqrt{D} \alpha^{D-1} \sqrt{1 - \alpha^2} \vb{e}_1^T \\
    \times & \identity & \sqrt{D}\alpha^{D-1} \sqrt{1 - \alpha^2} & A_S \\
    \times & \times & 1 & 0 \\
    \times & \times & \times & \identity_{n-1}
\end{bmatrix}
,\end{align*}
where $A_S = \mathrm{diag}(-\alpha^D + (D-1) \alpha^{D-2} (1 - \alpha^2), \alpha^{D-1},\dots,\alpha^{D-1})$.
Define the two matrices
$$
Z = \begin{bmatrix}
    \alpha^D & \sqrt{D} \alpha^{D-1} \sqrt{1 - \alpha^2} \vb{e}_1^T \\
    \sqrt{D}\alpha^{D-1} \sqrt{1 - \alpha^2} \vb{e}_1 & A_S \\
\end{bmatrix},\quad
Z^\prime := \begin{bmatrix}
    \alpha^D & \sqrt{D} \alpha^{D-1} \sqrt{1 - \alpha^2} \\
    \sqrt{D}\alpha^{D-1} \sqrt{1 - \alpha^2} & -D\alpha^D + (D-1)\alpha^{D-2} \\
\end{bmatrix}
.$$
The eigenvalues of $G_S$ are $1 \pm \sigma(Z)$.
Due to the sparse structure of $Z$, its singular values are $\alpha^{D-1}$ and the singular values of $Z'$. Since $Z'$ is symmetric, its eigenvalues and singular values coincide. We factor out $\alpha^{D-2}$ and compute the eigenvalues in terms of the trace $\tau$ and determinant $\Delta$. This gives $\tau = (D-1)(1-\alpha^2)$, $\Delta = -\alpha^2$, and
$
\lambda_1(Z') = \frac{\alpha^{D-2}}{2} \left( \tau + \sqrt{\tau^2 - 4\Delta} \right)$ and
$\lambda_2(Z') = \frac{\alpha^{D-2}}{2} \left( \tau - \sqrt{\tau^2 - 4\Delta} \right)
.$
Finally, we compare the eigenvalues of $G_S$ to the extreme eigenvalues of $G_\perp$.
Since $\alpha^2 \le 1$ and $D \ge 3$,
$$
4 \tau \ge 4(1 + \Delta)
\quad\Rightarrow\quad
\tau^2 - 4\Delta \ge \tau^2 - 4\tau + 4
\quad\Rightarrow\quad
\sqrt{\tau^2 - 4\Delta} \ge 2 - \tau
\quad\Rightarrow\quad
\frac{1}{2} (\tau + \sqrt{\tau^2 - \Delta}) \ge 1.
$$
Hence, $G_S$ has at least one eigenvalue less than or equal to the smallest of $G_\perp$, namely
$1 + \lambda_2(Z') \le 1 + \alpha^{D-2}$ if $\alpha^{D-2}$ is negative, and 1 - $\lambda_1(Z') \le 1 - \alpha^{D-2}$ otherwise. This shows that the smallest singular value of~$\widetilde{T}^{\mathcal{S}}$ in \cref{eq:TerraciniSegreBlocks} is a singular value of~$T_{\mathpzc{A}_1,\mathpzc{A}_2}^{\mathcal{V}}$, as required.
\end{proof}

\subsection{Numerical experiments}
\label{sec:experiments}
We tested \cref{conjecture:WaringCondition} for third order tensors. For $n = 3 \dots 18$, we generated 500 random symmetric rank $R$ decompositions $\sum_{r=1}^R \vb{a}_r^{\otimes 3}$ where $\vb{a}_r \sim \mathcal{N}(0,\identity_n)$ using Julia v1.6 \cite{Bezanson2017}. For each decomposition, we computed the two condition numbers. By dimensionality arguments, the condition number can only be finite if
$
Rn <
{n + 2 \choose 3}
$, where the right-hand side is the dimension of the space of symmetric $n \times n \times n$ tensors \cite{Landsberg2012}. We tested all values of $R$ below this upper bound.

\begin{figure}
    \centering
    \includegraphics[width=0.7\textwidth]{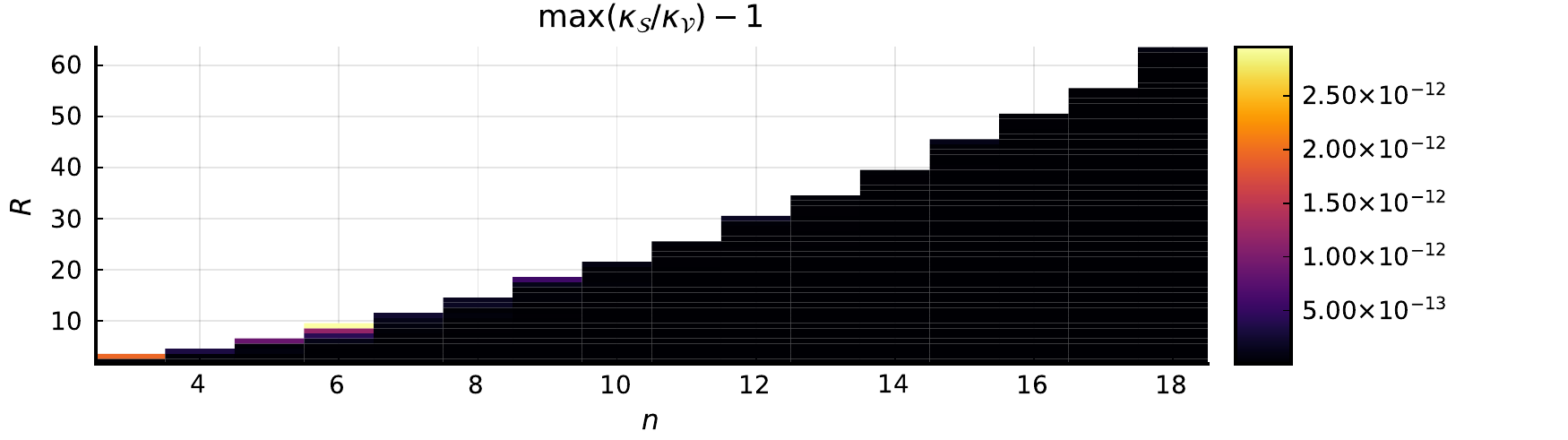}
    \caption{Ratio between the condition numbers of the CPD and WD of an $n \times n \times n$ symmetric tensor of rank $R$. The displayed value is the maximum over 500 test cases.}
    \label{fig:condratio_max}
\end{figure}

\Cref{fig:condratio_max} shows the ratio between the condition number of the CPD and the WD. A priori, it can never be less than 1. In practice, numerical computations would sometimes find a ratio of $1 - 10^{-11}$ or less. This suggests that ratios exceeding 1 by less than $10^{-11}$ can be explained by numerical roundoff. All measurements lie below this threshold.

\bibliographystyle{elsarticle-num}
\bibliography{library}

\clearpage
\appendix

\section{The condition number of the partially symmetric decomposition}
\label{sec:partiallySymmetric}

In this section, we present a generalisation of \cref{prop:WaringCondInvariance} to the partially symmetric case. We say that a tensor $\mathpzc{A}$ is partially symmetric if it is invariant under the permutation of some (but not all) of its indices. When this symmetry constraint is imposed on the summands in its CPD, the CPD is called a \textit{partially symmetric rank decomposition (PSRD)}.
Write the size and degree of the tensors as $\vb{n} = (n_1,\dots,n_K)$ and $\vb{d} = (d_1,\dots,d_K)$, respectively. Then partially symmetric tensors of rank 1 form the image of the map
\begin{align*}
\Phi: \mathbb{R} \setminus \{0\} \times \mathbb{S}^{n_1 - 1} \times \dots \times \mathbb{S}^{n_K - 1} & \rightarrow \mathbb{R}^{n_1 \times \dots \times n_1 \times \dots \times n_K} \\
(\alpha, \vb{a}_1,\dots,\vb{a}_K)
& \mapsto \alpha \vb{a}_1^{\otimes d_1} \otimes\dots\otimes \vb{a}_K^{d_K}
.\end{align*}
The image of $\Phi$ is known as the Segre-Veronese manifold $\mathcal{SV}_{\vb{n}, \vb{d}}$ \cite{Landsberg2012}. Analogous to the $Q$-WD, a $(Q_1,\dots,Q_K)$-PSRD is a PSRD of the form $\mathpzc{A} = \sum_{r=1}^R (Q_1^{\otimes d_1} \otimes \dots \otimes Q_K^{\otimes d_K}) \mathpzc{G}_r$ where each $Q_k$ has orthonormal columns and $\mathpzc{G}_r \in \mathcal{SV}_{\vb{m}, \vb{d}}$ where $\vb{m} < \vb{n}$ elementwise. We write $\mathcal{W} = (Q_1^{\otimes d_1} \otimes \dots \otimes Q_K^{\otimes d_K})(\mathcal{SV}_{\vb{m}, \vb{d}})$.

To determine the condition number, we apply \cref{eq:GeneralTerraciniMtx} to the PSRD. The derivative of $\Phi$ at any point $\mathpzc{A}_r = \Phi(\alpha, \vb{a}_1, \dots, \vb{a}_K)$ is
\begin{align*}
\mathrm{d} \Phi(\dot{\alpha}, \dot{\vb{a}}_1, \dots, \dot{\vb{a}}_K)
& =
\dot{\alpha}_r \bigotimes_{k=1}^K \vb{a}_k^{\otimes d_k}
+ \sum_{k=1}^K
\left(
\sum_{d=1}^{d_k}
    \dot{\vb{a}}_k \otimes_{d} \vb{a}_k^{\otimes d_k - 1}
\right)
\otimes_k
\left(
    \bigotimes_{k'\ne k} \vb{a}_{k'}^{\otimes d_{k'}}
\right)
.\end{align*}
If $U(\vb{a}_k)$ spans an orthonormal basis of $T_{\vb{a}_k} \mathbb{S}^{n_k - 1}$, the tangent space to $\mathcal{SV}_{\vb{n},\vb{d}}$ is the column space of
\begin{align}
    \label{eq:TerraciniPartSym}
T_{\mathpzc{A}_r}^{\mathcal{SV}_{\vb{n},\vb{d}}} &:=
\left[
    \bigotimes_{k=1}^K \vb{a}_k^{\otimes d_k}
\quad
\left[
    \frac{1}{\sqrt{d_k}}
    \left(\sum_{d=1}^{d_k} U(\vb{a}_k) \otimes_{d} \vb{a}_k^{\otimes d_k - 1}\right)
\otimes_k \left(
    \bigotimes_{k'\ne k} \vb{a}_{k'}^{\otimes d_{k'}}
\right)
\right]_{k=1}^K
\right]
.\end{align}
Observe that all $K+1$ blocks of this matrix have orthonormal columns and are pairwise orthogonal by construction of $U_k$. Therefore, the condition number of any PSRD can be computed using \cref{eq:GeneralTerraciniMtx} where the blocks in the Terracini matrix are as in \cref{eq:TerraciniPartSym}.
Now we can present a generalisation of \cref{prop:WaringCondInvariance}.

\begin{theorem}
    \label{prop:PSRDcondInvariance}
    Let $\mathpzc{G} = \mathpzc{G}_1 + \dots + \mathpzc{G}_R$ be a PSRD with summands in $\mathcal{SV}_{\vb{m}, \vb{d}}$. For $k = 1,\dots,K$, take $Q_k \in \mathbb{R}^{n_k \times m_k}$ with orthonormal columns and set $\mathpzc{A}_r := (Q_1^{\otimes d_1} \otimes \dots \otimes Q_K^{\otimes d_K}) \mathpzc{G}_r$. Then
    $$
    \kappa_{\mathcal{SV}_{\vb{n}, \vb{d}}}(\mathpzc{A}_1,\dots,\mathpzc{A}_R) \le \sqrt{\max \vb{d}}\,
    \kappa_{\mathcal{SV}_{\vb{m}, \vb{d}}}(\mathpzc{G}_1,\dots,\mathpzc{G}_R)
    .$$
    Similarly, for $k = 1,\dots,K$, let $U_k \in \mathbb{R}^{\ell_k \times m_k}$ have orthonormal columns and $\mathpzc{B}_r := (U_1^{\otimes d_k} \otimes\dots\otimes U_k^{\otimes d_k}) \mathpzc{A}_r$ for $r=1,\dots,R$. If $\min(\ell_k, n_k) > m_k$ for all $k$, then
    $$
    \kappa_{\mathcal{SV}_{\vb{n},\vb{d}}}(\mathpzc{A}_1,\dots,\mathpzc{A}_R)
    =
    \kappa_{\mathcal{SV}_{\widetilde{\vb{n}},\vb{d}}}(\mathpzc{B}_1,\dots,\mathpzc{B}_R)
    .$$
\end{theorem}
\begin{remark}
    The case $K = 1$ is exactly \cref{prop:WaringCondInvariance}. The case $d_1 = \dots = d_K = 1$ is a statement about the CPD. In this case the theorem reads $\kappa_{\mathcal{S}_{\vb{n}, K}}(\mathpzc{A}_1,\dots,\mathpzc{A}_R) = \kappa_{\mathcal{S}_{\vb{n}, K}}(\mathpzc{G}_1,\dots,\mathpzc{G}_R)$, which is a special case of \cite[Theorem 5.1]{Dewaele2021}.
\end{remark}

\begin{proof}
    For each $r$, let $\mathpzc{G}_r = \alpha_r (\vb{g}_1^r)^{\otimes d_1} \otimes\dots\otimes (\vb{g}_K^r)^{\otimes d_K}$ with $\alpha_r \ne 0$ and $\vb{g}_k^r \in \mathbb{S}^{m_k - 1}$ for all $k$. Let $\vb{a}_k^r = Q_k \vb{g}^r_k$ and define $U^r_k$ so that $[\vb{g}^k_r \quad U^k_r] \in \mathbb{R}^{m_k \times m_k}$ is orthogonal. Construct $T_{\mathpzc{G}_r}^{\mathcal{SV}_{\vb{n},\vb{d}}}$ by applying \cref{eq:TerraciniPartSym} to $\mathpzc{G}_r$.
    Complete each $Q_k$ to an orthonormal basis $[Q_k \quad Q_k^\perp]$ of $\mathbb{R}^{n_k}$. If $n_k = m_k$, $Q_k^\perp$ is an $n_k \times 0$ matrix. The columns of $U(\vb{a}_k^r) := [Q_k U_k^r\quad Q^\perp_k]$ form an orthonormal basis of $T_{\vb{a}^r_k} \mathbb{S}^{n_k - 1}$. For each $r$, these can be substituted into \cref{eq:TerraciniPartSym} applied to $\mathpzc{A}_1,\dots,\mathpzc{A}_R$, respectively. Similarly to the symmetric case, this gives
    $$
    T_{\mathpzc{A}_r}^{\mathcal{SV}_{\vb{n},\vb{d}}} = \left[
        T_r \quad T_r^{1\perp} \quad \cdots \quad T_r^{K \perp}
    \right]
    \quad\text{where}\quad
    T_r = \left( \bigotimes_{k=1}^K Q_k^{\otimes d_k} \right)
    T_{\mathpzc{G}_r}^{\mathcal{SV}_{\vb{m},\vb{d}}}
    \quad\text{and}\quad
    $$$$
    T_r^{k\perp} =
        \frac{1}{\sqrt{d_k}} \left(\sum_{d=1}^{d_k} Q_k^\perp \otimes_d (\vb{a}^r_k)^{\otimes d_k - 1}\right)
    \otimes_k \left(
        \bigotimes_{k'\ne k} (\vb{a}^r_{k'})^{\otimes d_{k'}}
    \right)
    $$
    for each $r$ and $k$. Define $T = [T_r]_{r=1}^R$ and $T^{k\perp} = [T_r^{k\perp}]_{r=1}^R$. Observe that these $K+1$ matrices are pairwise orthogonal since $Q_k^T Q_k = 0$ and $\vb{a}_k^r \in \mathrm{span}\,Q_k$. Furthermore, note that $T_{\mathpzc{A}_r}^{\mathcal{SV}_{\vb{n},\vb{d}}} = [T \quad T^{1\perp} \quad\cdots\quad T^{K\perp}]$ up to a column permutation. Finally, $T = \left( \bigotimes_{k=1}^K Q_k \right) T_{\mathpzc{G}_1,\dots,\mathpzc{G}_R}^{\mathcal{SV}_{\vb{m},\vb{d}}}$ has the same $T_{\mathpzc{G}_1,\dots,\mathpzc{G}_R}^{\mathcal{SV}_{\vb{m},\vb{d}}}$ by orthogonality. The combination of these three observations implies that the singular values of $T_{\mathpzc{A}_r}^{\mathcal{SV}_{\vb{n},\vb{d}}}$ are the union of the singular values of $T,\, T^{1\perp},\dots,T^{K\perp}$ separately.
    Consequently, it suffices to show that for each $k$ we have $\sigma_{\min}(T^{k\perp}) \ge \sigma_{\min}(T_{\mathpzc{G}_1,\dots,\mathpzc{G}_R}^{\mathcal{SV}_{\vb{m},\vb{d}}}) / \sqrt{d_k}$.
    To do this, we compute the Gramian $(T^{k\perp})^T T^{k\perp}$.
    Define the following auxiliary matrices:
    $$
    S^k_r = \frac{1}{\sqrt{d_k}} \left(\sum_{d=1}^{d_k} Q_k^\perp \otimes_d (\vb{a}^r_k)^{\otimes d_k - 1}\right), \qquad
    A^k_r = \bigotimes_{k'\ne k} (\vb{a}^r_{k'})^{\otimes d_{k'}}, \qquad
    G^k_r = \bigotimes_{k'\ne k} (\vb{g}^r_{k'})^{\otimes d_{k'}}
    .$$
    This allows us to write $T^{k\perp} = S_r^k \otimes_k A_r^k$.
    For general $r_1,r_2$, the inner products between the columns of $S_{r_1}^k$ and $S_{r_2}^k$ are
    $$(S_{r_1}^k)^T S_{r_2}^k =
    \left\langle
    \vb{a}^{r_1}_k,\,
    \vb{a}^{r_2}_k
    \right\rangle^{d_k - 1} \identity_{n_k - m_k}
    =
    \left\langle
        \vb{g}^{r_1}_k
    ,\,
        \vb{g}^{r_2}_k
    \right\rangle^{d_k - 1} \identity_{n_k - m_k}
    .$$
    Hence, if we replace the factors $S_r^k$ in $T_r^{k\perp}$ by $\identity_{n_k - m_k} \otimes (\vb{g}_k^r)^{\otimes d_k - 1}$, the Gramian remains unchanged. Similarly, $(A_{r_1}^k)^T A_{r_2}^k = (G_{r_1}^k)^T G_{r_2}^k$ for all $r_1$ and $r_2$, so that we can replace each $A_r^k$ in $T^{k\perp}$ by $G_r^k$. Define
    $$
    \hat{T}^{k\perp} := \left[
        \identity_{n_k - m_k} \otimes (\vb{g}_k^r)^{\otimes d_k - 1} \otimes_k G_k^r
    \right]_{r=1}^R
    \quad\text{and}\quad
    \widetilde{T}^{k\perp} := \left[
        [\vb{g}_k^r \quad U_k^r] \otimes (\vb{g}_k^r)^{\otimes d_k - 1} \otimes_k G_k^r
    \right]_{r=1}^R
    .$$
    $\hat{T}^{k\perp}$ is $T^{k\perp}$ with the aforementioned replacements applied. Since $[\vb{g}_k^r \quad U_k^r]$ is orthogonal, the singular values of $\hat{T}^{k\perp}$ and $\widetilde{T}^{k\perp}$ are the same up to multiplicities \cite[Lemma 5.3]{Dewaele2021}. Hence, for the purpose of comparing singular values, we can proceed with $\widetilde{T}^{k\perp}$ instead of $T^{k\perp}$.

    Next, we also modify $T_{\mathpzc{G}_1,\dots,\mathpzc{G}_R}^{\mathcal{SV}_{\vb{m}, \vb{d}}}$.
    First, take the following subset of its columns:
    $$
    T^k
    :=
    \left[
        \bigotimes_{k=1}^K (\vb{g}_k^r)^{\otimes d_k}
        \quad
        \frac{1}{\sqrt{d}_k}
        \left( \sum_{d=1}^{d_k}
            U_k^{r\perp} \otimes_d (\vb{g}_k^r)^{\otimes d_k - 1}
        \right)
        \otimes_k G_r^k
    \right]_{r=1}^R.
    $$
    The first column of the $r$th block is $\bigotimes_{k=1}^K (\vb{g}_k^r)^{\otimes d_k} = (\vb{g}_k^r)^{\otimes d_k} \otimes_k G_r^k$. Define $\widetilde{T}^k$ as a modification of $T^k$ where these $R$ columns are scaled up by $\sqrt{d_k}$. Rearranging the columns gives
    $$
    \widetilde{T}^k
    =
    \left[
        \frac{1}{\sqrt{d_k}} \left(
            \sum_{d=1}^{d_k}
            [\vb{g}_k^r\quad U_k^{r\perp}] \otimes_d (\vb{g}_k^r)^{\otimes d_k - 1}
        \right)
        \otimes_k G_r^k
    \right]_{r=1}^R
    .$$
    Since $T^k$ is a submatrix of $T_{\mathpzc{G}_1,\dots,\mathpzc{G}_R}^{\mathcal{SV}_{\vb{m}, \vb{d}}}$, we have $\sigma_{\min}(T_{\mathpzc{G}_1,\dots,\mathpzc{G}_R}^{\mathcal{SV}_{\vb{m}, \vb{d}}}) \le \sigma_{\min}(T^k)$. Because of how we defined $\widetilde{T}^k$, we also have $\sigma_{\min}(T^k) \le \sigma_{\min}(\widetilde{T}^k)$.
    From here on, we can compare the singular values of $\widetilde{T}^k$ and $\widetilde{T}^{k\perp}$
    the same way as their counterparts in the proof of \cref{prop:WaringCondInvariance}. This completes the proof.
\end{proof}

\end{document}